\documentclass[11pt]{amsart}

\catcode`\@=11
\usepackage{amsmath,amsfonts,amscd,amssymb}

\theoremstyle{plain}

\newtheorem{theorem}[equation]{Theorem}
\newtheorem{proposition}[equation]{Proposition}
\newtheorem{lemma}[equation]{Lemma}
\newtheorem{corollary}[equation]{Corollary}
\theoremstyle{definition}
\newtheorem*{acknowledgements}{Acknowledgements}
\newtheorem{remark}[equation]{Remark}
\newtheorem{definition}[equation]{Definition}
\newtheorem{notation}[equation]{Notation}
\newtheorem*{notation*}{Notation}
\newtheorem{example}[equation]{Example}
\newtheorem*{example*}{$\bullet$}

\def\legendre#1#2{\Bigl(\frac{#1}{#2}\Bigr)}
\def\tlegendre#1#2{\bigl(\frac{#1}{#2}\bigr)}

\newcommand{\vabove}[2]{\genfrac{}{}{0pt}{}{#1}{#2}}

\def\leftchoice#1#2#3#4{{\def\arraystretch{0.7}
\Bigl\{\!\!\begin{array}{ll}
   \scriptstyle #1,\!\!\!&\scriptstyle #2\cr
   \scriptstyle #3,\!\!\!&\scriptstyle #4\end{array}}}
\def\bigleftchoice#1#2#3#4{{\def\arraystretch{1}
   \biggl\{\!\!\begin{array}{ll}#1,&#2\cr#3,&#4\end{array}}}

\def\smz{{\scriptstyle z}}
\def\sms{{\scriptstyle \ast}}
\def\smm#1{{\scriptstyle \raise0.17em\hbox to 0.2em{\hrulefill}\hbox to 0.1em{\hfill}#1}}

\let\iff\Leftrightarrow
\let\liff\Longleftrightarrow
\def\injects{\lhook\joinrel\rightarrow}

\let\iso\cong

\def\tinymatrix#1{\begingroup
  \font\txtfnt=cmr7\textfont0=\txtfnt
  \baselineskip 8pt
  \def\arraystretch{0.5}
  \def\arraycolsep{2pt}
  \left(\!
  \begin{array}{ccccccccc}
  #1
  \end{array}
  \!\right)
  \endgroup}

\def\beq{$$\begin{array}{llllllllllllllll}}
\def\eeq{\end{array}$$}

\def\comment{}
\def\endcomment{}

\def\<{\raise0.5pt\hbox{$\,\scriptstyle<\,$}}

\def\bb@symb#1|#2{\expandafter\def\csname #2#1\endcsname{\mathbb{#1}}}
\def\bbsymbols#1#2{\@for\@tmpz:=#2\do{\expandafter\bb@symb\@tmpz|{#1}}}

\def\cal@symb#1|#2{\expandafter\def\csname #2#1\endcsname{\mathcal{#1}}}
\def\calsymbols#1#2{\@for\@tmpz:=#2\do{\expandafter\cal@symb\@tmpz|{#1}}}

\def\dmth@p#1|{\expandafter\let\csname#1\endcsname\relax
  \expandafter\DeclareMathOperator\csname#1\endcsname{#1}}
\def\operators#1{\@for\@tmpz:=#1\do{\expandafter\dmth@p\@tmpz|}}

\operators{Tr,Gal,Hom,coker,Res,Frob,GL,PGL,SL,Aut,Sy,id}
\DeclareMathOperator\vchar{char}
\let\tr\Tr

\bbsymbols{}{R,C,Z,N,Q,F}
\calsymbols{}{O}
\catcode`\@=12

\theoremstyle{plain}
\newtheorem{problem}[equation]{Problem}
\calsymbols{c}{A,B,C}
\bbsymbols{}{A}
\operators{Cent,PSL,PGL,PGSp,GSp,Disc,Sp}
\def\Cy{{\rm C}}
\def\Ma{{\rm M}}
\def\Di{{\rm D}}
\def\Sy{{\rm S}}
\def\Alt{{\rm A}}

\csname@addtoreset\endcsname{equation}{section}
\begin{document}

\def\a{{\boldsymbol a}}
\def\b{{\boldsymbol b}}
\def\p{{\mathfrak p}}
\def\q{{\mathfrak q}}
\def\e#1#2{e_{#1}^{\scriptscriptstyle #2}}
\def\eFa{\e{\a}{F}}
\def\m{{{\Gamma\hskip-0.8pt}}}
\def\mm#1#2#3{\m_{#1,#2}^{\scriptscriptstyle #3}}
\def\mmaD{\mm{\a}{D}{F}}
\def\cM#1#2#3{{M}^{#1}_{#2,#3}}

\comment

\title{Identifying Frobenius elements in Galois groups}
\author{Tim and Vladimir Dokchitser}
\address{Robinson College, Cambridge CB3 9AN, United Kingdom}
\email{t.dokchitser@dpmms.cam.ac.uk}
\address{Emmanuel College, Cambridge CB2 3AP, United Kingdom}
\email{v.dokchitser@dpmms.cam.ac.uk}
\llap{.\hskip10cm}\vskip-4mm

\begin{abstract}
We present a method to determine Frobenius elements
in arbitrary Galois extensions of global fields,
which may be seen as a generalisation of Euler's criterion.
It is a part of the general question how to
compare splitting fields
and identify conjugacy classes in Galois groups,
that we will discuss as well.
\end{abstract}

\date{April 20, 2011}

\maketitle

\vskip-8mm\llap{.\hskip10cm}

\begingroup
\smaller[1]
\tableofcontents
\endgroup

\vskip-10mm\llap{.\hskip10cm}

\section{Introduction}
\label{sintro}

Take a Galois extension $L/\Q$.
Associated to each (unramified) prime $p$ is a Frobenius element $\Frob_p$,
an element of the Galois group that reduces to
$x\mapsto x^p$ modulo a prime above $p$.
In the setting when $L$ is the splitting field of a polynomial $f$,
this element is intimately connected to the factorisation of $f$ mod $p$:
viewed as a permutation of the roots,
$\Frob_p$ is a product of disjoint cycles whose lengths are the degrees
of the irreducible factors.

In this paper, we address the question how to determine $\Frob_p$.
Generally, we study the problem how to compare splitting fields
and identify conjugacy classes in Galois groups
(see \S\ref{sordering}-\ref{sts}).
Our motivation was computing L-series of Artin
representations for arbitrary Galois groups, which require
the knowledge of Frobenius elements
at all primes (see Remark~\ref{artinl} and Example~\ref{pgsp43}).
Obtaining them directly from the definition is impractical
unless $L$ either has small degree or
is particularly simple to work with.

Let us briefly illustrate the various standard techniques for computing
Frobenius elements.
As before, $L$ is the splitting field of a polynomial~\hbox{$f\!\in\!\Z[x]$},
and we write $G=\Gal(L/\Q)$.

\begin{example*}[Quadratic fields] Suppose $f(x)=x^2-d$, so $L=\Q(\sqrt d)$.
For a prime $p\nmid 2d$, the Frobenius element is given by the Legendre
symbol:
$$
 \Frob_p = \id  \quad\iff\quad
  f(x)\bmod p \text{ \ reducible}
  \quad\iff\quad
  \legendre dp=1.
$$
There are two essentially different methods to compute it:
\begin{itemize}
\item[(A)] Euler's criterion $\tlegendre dp=(-1)^{\frac{p-1}2}$.
\item[(B)] Quadratic reciprocity.
\end{itemize}
\end{example*}

\begin{example*}[Kummer extensions]
Suppose $f(x)\!=\!x^3-2$, so $L\!=\!\Q(\zeta_3,\sqrt[3]2)$ and~$G\!=\!\Sy_3$.
For $p\ne 2,3$ the number of cube roots of $2$ mod $p$ determines whether
$\Frob_p$ is $\id$, a 3-cycle or a transposition.
It it easy to see that the last case is equivalent to $p\equiv 2\mod 3$.
There are analogues of both (A) and (B) to distinguish between the
first two cases:

\smallskip\noindent
(A) Euler's criterion: since $\F_p^\times$ is cyclic,

\noindent
\begin{tabular}{lllllll}
&&2 is a cube  mod $p$ &$\iff$&  $2^{\frac{p-1}3}\equiv1\mod p$.\cr
&&2 not a cube mod $p$ &$\iff$&  $2^{\frac{p-1}3}$ is another third
                  root of unity $z\in\F_p$.\\[2pt]
\end{tabular}

\noindent
To link this criterion to our main theorem below, let us rephrase it:
let $M=\tinymatrix{0&0&2\cr1&0&0\cr0&1&0\cr}\in\GL_3(\F_p)$, so that $M^3=2$. Then
\begingroup
\smaller[0]
$$
\begin{array}{llllllll}
\Frob_p=\id &\iff\!& M^{p-1}\!=\!\tinymatrix{1&0&0\cr0&1&0\cr0&0&1\cr} &\iff&
  \frac13\tr M^{p-1}=1,\cr
\Frob_p\in[(123)] &\iff\!& M^{p-1}\!=\!\tinymatrix{\smz&0&0\cr0&\smz&0\cr0&0&\smz\cr} &\iff&
  \frac13\tr M^{p-1}\text{ satisfies $t^2\!+\!t\!+\!1\!=\!0$},\cr
\Frob_p\in[(12)] &\iff\!& M^{p-1}\!=\!\tinymatrix{0&0&\sms\cr\sms&0&0\cr0&\sms&0\cr} &\iff&
  \frac13\tr M^{p-1}=0.\cr
\end{array}
$$
\endgroup

\noindent
(B) Class field theory over $\Q(\zeta_3)$:

Factorise $p=(a+b\zeta_3)(a+b\bar\zeta_3)$. Then 2 is a cube mod $p$
if and only if the ideal $(a+b\zeta_3)$ splits in $L$, and
class field theory says that this is a congruence condition on $a$ and $b$.
In fact, it is easy to verify that
\begin{center}
2 is a cube mod $p$ $\ \iff\ $ $a+b\zeta_3\equiv \pm1,\pm\zeta_3$ or $\pm\zeta_3^2$ mod 6.
\end{center}
\end{example*}

\begin{example*}[Modular forms, see \cite{Z123}, \S4.3]
Suppose $f(x)=x^3-x-1$, so $G=\Sy_3$ and $L$ is the Hilbert class field
of $\Q(\sqrt{-23})$. Let $\rho$ be the 2-dimensional irreducible
representation of $G$. It has an associated Artin L-series
$$
  L(\rho,s) = \sum_{n=1}^\infty \frac{a_n}{n^s},
$$
whose coefficient $a_p$ for a prime $p\ne 23$ is $2, -1$ or $0$ depending
on whether $\Frob_p$ is trivial, a 3-cycle or a transposition.
The theory of modular forms tells us that
$$
  \sum_{n=1}^\infty a_n q^n = q\prod_{n=1}^\infty (1-q^n)(1-q^{23n}),
$$
and is a cusp form of weight 1, level 23 and character
$\tlegendre{\cdot}{23}$. Moreover, for all integers $n$ not divisible by 23,
$$
  a_n = \frac12\Bigl(\#\bigl\{x,y\!\in\!\Z \bigm| n\!=\!x^2\!+\!xy\!+\!6y^2\bigr\}\>-\>
                \#\bigl\{x,y\!\in\!\Z \bigm| n\!=\!2x^2\!+\!xy\!+\!3y^2\bigr\}\Bigr).
$$
Let us remark that in an arbitrary Galois group $G$, the L-series of
the irreducible representations of $G$ also pin down the Frobenius elements.
The global Langlands conjecture predicts that, as in this example,
all such L-series come from automorphic forms.
This is a massive conjectural generalisation of `method (B)'. Moreover,
like quadratic reciprocity and class field theory, this approach gives
expressions for the L-series coefficients $a_n$ that do not depend on
$n$ being prime. This is crucial for theoretical applications such as
analytic continuation of L-functions. (Note, however, that
formulae such as the one above are not practical for numerically computing
Frobenius elements.)
\end{example*}

The purpose of this paper is to extend `method (A)' to
arbitrary Galois groups. Here is an illustration for cubic polynomials
of the type of criterion that we obtain.
Note its similarity to the Kummer case.

\begin{example*}[General cubic]
Suppose $f(x)=x^3+bx+c$.
Pick a prime $p\nmid 3b\Delta_f$, where $\Delta_f\!=\!-4b^3\!-\!27c^2$ is the discriminant
of $f$.
Let $M=\tinymatrix{0&0&\smm c\cr1&0&\smm b\cr0&1&0\cr}\in\GL_3(\F_p)$.~Then
$$
\begin{array}{lll}
  \text{$f(x)$ has 3 roots mod $p$}    &\iff& \tr M^{p+1}=-2b,\cr
  \text{$f(x)$ has 1 root mod $p$}     &\iff& \tr M^{p+1}
    \hbox{ satisfies $(t+2b)(t-b)^2=-\Delta$,}\cr
  \text{$f(x)$ is irreducible mod $p$} &\iff& \tr M^{p+1}=b.\cr
\end{array}
$$
This can be easily checked by hand; alternatively, see Example \ref{gencubus}.
\end{example*}

Our main result for Frobenius elements
is the following generalisation of Euler's criterion.
Note that taking the class of $x$ in $\F_q[x]/f(x)$ is the same as
taking a matrix $M$ with characteristic polynomial $f(x)$, like in the
examples above.

\begin{theorem}
\label{imain}
Let $K$ be a global field and $f(x)\in K[x]$ a separable
polynomial with Galois group $G$.
There is a polynomial $h(x)\in K[x]$ and polynomials $\m_C\in K[X]$
indexed by the conjugacy classes $C$ of $G$ such that
$$
  \Frob_\p\in C \quad\iff\quad
    \m_C\Bigl(\tr_{\frac{\F_q[x]}{f(x)}/\F_q}(h(x)x^q)\Bigr)=0\mod\p
$$
for almost all primes $\p$ of $K$;
here $\F_q$ is the residue field at $\p$.
\end{theorem}
This is proved in \S\ref{sfrob}, see Theorem \ref{frobmain}.
Usually one can take $h(x)=x^2$ (see below),
in particular $\Tr(x^{q+2})$ then determines the conjugacy class of $\Frob_\p$.
In \S\ref{sexaab} we explain
how the theorem recovers classical formulae for Frobenius elements
in cyclotomic and Kummer extensions.
In \S\ref{sexanab}
we give explicit examples for non-abelian Galois groups,
including general
cubics, quartics and quintics with Galois group $\Di_{10}$.

The polynomials $\m_C$ are explictly given~by
$$
  \m_C(X) \>=\> \prod_{\sigma\in C}\, \bigl(X-
    \mathop{{\raise-3.5pt\hbox{\larger[3]$\Sigma$}}}\limits_{j=1}^n h(a_j) \sigma(a_j)
           \bigr),
$$
where $a_1,...,a_n$ are the roots of $f$ in some splitting field.
The `almost all primes' in the theorem are those not dividing the denominators of
the coefficients of $f$, its leading coefficient and the resultants
$\Res(\m_C,\m_{C'})$ for $C\ne C'$; the latter
simply says that the $\m_C$ mod $\p$ are pairwise coprime.
(This condition always fails for ramified primes, see Remark~\ref{rampr}.)
Finally, the only constraint on the polynomial $h$
is that the resulting $\m_C$
are coprime over $K$.
This holds for almost all $h$, in the sense that
the admissible ones of degree at most $n-1$ form a Zariski dense
open subset of $K^n$. Also, a fixed $h$ with $1<\deg h<n$ (e.g. $h(x)=x^2$)
will work for almost all $f$ that define the same field
(see~\S\ref{sugly}).

\begin{remark}
The method of using polynomials in the roots of $f$
to recognise conjugacy classes
is also used in ``Serre's trick'' for alternating groups.
For example, $G=A_5$ has 5 conjugacy classes and all but the
two classes of 5-cycles have their own cycle type.
(Recall that the cycle type of Frobenius can be recovered from the
degrees of the factors of the defining quintic $f$ mod~$p$;
in practice, these are readily determined by computing
$\gcd(x^{p^d}\!-\!x,f(x))$ for $d=1,2$.)
It was pointed out by Serre (see Buhler \cite{Buh} p. 53)
that the classes of 5-cycles can be distinguished by evaluating
the square root of the discriminant of $f$ modulo $p$; cf. Example \ref{serre1}.
This has been generalised by Roberts \cite{Rob} to all alternating groups,
and was used for instance by Booker~\cite{Boo} in his work on
L-series for icosahedral representations.
\end{remark}

Finally, let us illustrate some of the ideas in our approach
to Frobenius elements with a simple case:

\begin{example}
\label{iexample}
The polynomial $f(x)=x^5+2x^4-3x^3+1$ has Galois group $G=\Di_{10}$ over $K=\Q$.
If we number its complex roots by
$$
  a_1\approx-3.01, \quad
  a_2\approx-0.35\!-\!0.53i,\quad
  a_3\approx 0.85\!-\!0.31i,\quad
  a_4=\overline{a_3},\quad
  a_5=\overline{a_2},
$$
then $G$ is generated by the 5-cycle $(12345)$ and complex conjugation
$(25)(34)$. It is easy to see that $f(x)$ is irreducible over $\F_2$,
so $\Frob_2\in G$ is in one of the two conjugacy classes of 5-cycles,
either $[(12345)]$ or $[(12345)^2]$. How can we check which one it is?

Consider the expressions,
\beq
  n_1 &=& a_1 a_2 + a_2 a_3 + a_3 a_4 + a_4 a_5 + a_5 a_1,\cr
  n_2 &=& a_1 a_3 + a_2 a_4 + a_3 a_5 + a_4 a_1 + a_5 a_2.
\eeq
If we think of $G$ as the group of symmetries of a pentagon, the sums
are taken over all edges and over all diagonals respectively.
Therefore they are clearly $G$-invariant, i.e. rational
numbers. Moreover, as $a_i$ are algebraic integers, $n_1$ and $n_2$
are in fact integers, readily recognised
from their complex approximations as being 2 and~$-5$.

Now suppose $b_1$ is a root of $f(x)$ in $\F_{2^5}$, and
$b_i=b_{i-1}^2$ for $i=2,3,4,5$ are its other roots ordered by the action
of the Frobenius automorphism. Then
$$
  N = b_1 b_2 + b_2 b_3 + b_3 b_4 + b_4 b_5 + b_5 b_1
$$
is in $\F_2$. By considering the reduction modulo a prime $\q$ above 2
in the splitting field, we see that if $\Frob_\q$ is $(12345)$ or $(12345)^{-1}$,
then $n_1\equiv N\mod 2$. Similarly, if $\Frob_\q$ is $(12345)^2$ or $(12345)^3$,
then $n_2\equiv N\mod 2$. Computing in $\F_2^5$
(or noting that $N = \Tr_{\F_2[x]/f(x)}(x^3)$) we find that $N=0$,
so $\Frob_2$ must be in $[(12345)]$.

In the language of Theorem \ref{imain}, we took $h(x)=x$ and proved that
$$
  \m_{[(12345)]}=(X-2)^2 \qquad\text{and}\qquad  \m_{[(12345)^2]}=(X+5)^2
$$
distinguish between the two conjugacy classes of 5-cycles:
if $f(x)$ is irreducible mod $p$ (and $p\ne 7$, so that $2\not\equiv-5$), then
$$
  \Frob_p\in C\quad\iff\quad
  \m_C(\Tr_{\F_p[x]/f(x)}(x^{p+1}))=0\mod p.
$$
This choice of $h(x)$ was in some sense deceptively simple, because
the roots $n_i$ of the $\m_C$'s were integers.
(We used that the conjugacy classes of 5-cycles are self-inverse in
$\Di_{10}$.) 
Generally, these roots would be algebraic integers of degree $|C|$.
For example, $h(x)=x^2$ leads to
$$
  \m_{[(12345)]}=X^2+5X+18 \qquad\text{and}\qquad  \m_{[(12345)^2]}=X^2-11X+42,
$$
and $\Tr(x^{p+2})$ is a root of one of them whenever $f(x)\mod p$
is irreducible.
\end{example}

\begin{notation*}
Throughout the paper we use the following notation:

\par
\vskip 1mm
\noindent
\hskip-1mm\begin{tabular}{lll}
&$K$                & ground field\cr
&$f(x)$             & separable polynomial in $K[x]$ of degree $n$\cr
&$L$                & some extension of $K$ where $f$ splits completely\cr
&$\a=[a_1,...,a_n]$ & ordered roots of $f$ in $L$\cr
&$K(\a)$            & field generated by the $a_i$ over $K$ (a splitting field of $f$)\cr
&$G_\a$             & Galois group of $f$, considered as a subgroup of $\Sy_n$ \cr
&                   & via its permutation action on $[a_1,...,a_n]$.\cr
&$[\Psi]$           & conjugacy class of $\Psi\in G_\a$.\cr
&$\p$               & prime of $K$, when $K$ is a global field\cr
&$\F_q$             & residue field at $\p$\cr
&$\Frob_\p$         & any (arithmetic) Frobenius element at $\p$ in $G_\a$   \cr
&$\eFa, \m, \cM F\a\Psi$& see Definitions \ref{defe}, \ref{defm}, \ref{defcM}
                      and \ref{gammaC}. \\[2pt]
\end{tabular}

\noindent
Recall that a global field is a finite extension of either $\Q$ or $\F_p(T)$.
The Fro\-benius element in $\Gal(L/K)$ at $\p$ is
characterised by $\Frob_\p(x)\equiv x^q\mod\q$ for all $x\in L$
that are integral at some fixed prime $\q$ of $L$ above $\p$.
The element $\Frob_\p$ is well-defined
modulo inertia and up to conjugation. In particular, its conjugacy
class is well-defined if $\p$ is unramified in $L/K$.

The symmetric group $\Sy_n$ acts on $n$-tuples by
$$
  [c_1,...,c_n]^\sigma = [c_{\sigma^{-1}(1)},...,c_{\sigma^{-1}(n)}].
$$
It acts on the ring of polynomials in $n$ variables $K[x_1,...,x_n]$ by
$\sigma(x_i)=x_{\sigma(i)}$; thus, for a polynomial $F\in K[x_1,...,x_n]$,
$$
  F^\sigma([c_1,...,c_n]) = F([c_1,...,c_n]^{\sigma^{-1}}), 
$$
where $F([...])$ is the evaluation of $F$ on the $n$-tuple.
\end{notation*}

\begin{acknowledgements}
The first author is supported by a Royal Society
University Research Fellowship.
The second author would like to thank Gon\-ville \& Caius College, Cambridge,
where part of this research was carried~out.
\end{acknowledgements}

\section{Isomorphisms of splitting fields}
\label{sordering}

In this section we introduce our main tools. The reader who is only
interested in applications to Frobenius elements may skip to
\S\ref{sfrob} and prove Theorem \ref{frobmain} directly
(at the expense of not seeing the origins of $\m_C$).

As a motivation, consider the following general question:

\begin{problem}
\label{problem1}
Suppose we are given a
separable polynomial $f(x)\in K[x]$ of degree $n$ which splits completely
in $L\supset K$ and $L'\supset K$.
Given the roots $a_1,...,a_n$ and $b_1,...,b_n$
of $f$ in $L$ and $L'$, find a bijection between them that
comes from an isomorphism of splitting fields of $f$ inside $L$ and $L'$.
\end{problem}

\noindent
We assume that we know the Galois group of $f$ over $K$ as a permutation
group on the roots in $L$, but we do not want to construct the splitting fields
explicitly.
Instead, we will evaluate polynomials in $K[x_1,...,x_n]$ on the roots
in $L$ and $L'$ taken in various orders and try to extract information
out of the values (as in Example \ref{iexample}).

\begin{definition}
\label{defe}
For $F\!\in\! K[x_1,...,x_n]$ define the {\em evaluation map}
$
  \Sy_n\!\to\!K(\a)
$
by
$$  \eFa(\sigma) = F([a_1,...,a_n]^\sigma).
$$
\end{definition}

\begin{definition} 
For a subgroup $T$ of $\Sy_n$ a {\em $T$-invariant} $F$ is an element
of $K[x_1,...,x_n]$ whose stabiliser is precisely $T$.
\end{definition}

\begin{remark}
\label{invariantexists}
Any $F\in K[x_1,...,x_n]$ is evidently $T$-invariant if we take for
$T$ its stabiliser in $\Sy_n$. Also, any subgroup $T\< \Sy_n$ has a $T$-invariant, e.g.
$$
  F = \sum_{t\in T} m^t, \qquad m=x_1^{n-1}x_2^{n-2}\cdots x_{n-1},
$$
since clearly the stabiliser of $m$ in $\Sy_n$ is $\{1\}$.
\end{remark}

\begin{lemma}
\label{eFprop}
Let $F$ be a $T$-invariant and $\sigma,\tau\in \Sy_n$.
\begin{enumerate}
\item $\e{\a^\tau}{F}(\sigma)=\eFa(\sigma\tau)$ .
\item $g(\eFa(\sigma))=\eFa(\sigma g^{-1})$ for $g\in G_\a$.
\item The map $\eFa: \Sy_n\to K(\a)$ is constant on the right cosets $T\sigma$.
\end{enumerate}
\end{lemma}

\begin{proof}
(1)
$
  \e{\a^\tau}{F}(\sigma) =
  F((\a^\tau)^\sigma)=F(\a^{\sigma\tau})=\eFa(\sigma\tau).
$

\noindent
(2) For $g\in G_\a$,
\beq
  g(\eFa(\sigma)) \!\!\!&=\!\!\!&
  g(F([a_1,...,a_n]^\sigma)) \!\!\!&=\!\!\!& F([g(a_1),...,g(a_n)]^\sigma)\cr
    \!\!\!&=\!\!\!& F(([a_1,...,a_n]^{g^{-1}})^\sigma)
    \!\!\!&=\!\!\!& F([a_1,...,a_n]^{\sigma g^{-1}}) \quad=\>\, \eFa(\sigma g^{-1}).
\eeq
(3) For $\tau\in T$,
\beq
  \eFa(\tau\sigma) \!\!\!&=\!\!\!&
  F([a_1,...,a_n]^{\tau\sigma}) \!\!\!&=\!\!\!& F(([a_1,...,a_n]^\sigma)^\tau)
 \cr
    \!\!\!&=\!\!\!& F^{\tau^{-1}}([a_1,...,a_n]^\sigma)
    \!\!\!&=\!\!\!& F([a_1,...,a_n]^\sigma)  \quad=\>\,  \eFa(\sigma).
\eeq
\end{proof}

\begin{remark}
Part (3) of the lemma says that
the values of $F$ on the various permutations $\a^\sigma$ of the roots
are essentially the right cosets of $T$ in $\Sy_n$.
It may accidentally happen that the same value occurs on two right cosets,
but it is always possible to adjust the original polynomial $f$
to prevent this (see Lemma \ref{ugly}c).
Part (2) of Lemma \ref{eFprop} says that
the action of the Galois group $\Gal(K(\a)/K)$ on these values translates
into right multiplication by~$G_\a$. This motivates the following
\end{remark}

\begin{definition}
\label{defm}
For a double coset $D=T\sigma_0 G_\a$ in $\Sy_n$,
define the corresponding `minimal
polynomial'
$$
  \mm{\a}{\sigma_0}{F}=
  \mmaD(X) = \prod_{\sigma\in T\backslash D} (X-\eFa(\sigma)) \in K[X].
$$
By Lemma \ref{eFprop} (3), this is well-defined. 
\end{definition}

\begin{remark}
\label{reminj}
Note that by Lemma \ref{eFprop} (2),
$G_\a$ permutes the linear factors of $\mmaD$ transitively, so it is a power
of an irreducible polynomial in $K[X]$. If
$\eFa: T\,\backslash\,\Sy_n\to K(\a)$ is injective,
then $\mmaD(X)$ is irreducible, and hence the minimal polynomial
of $\eFa(\sigma_0)$.
\end{remark}

\begin{remark}
The point is that the $\mmaD(X)$ are $K$-rational objects, and they can
be used to compare different splitting fields:
\end{remark}

\begin{proposition}
Let $\a, \b$ be orderings of roots of $f$ in two splitting fields of~$f$,
and let $\phi: K(\a)\to K(\b)$ be an isomorphism.
If $\eFa: T\,\backslash\,\Sy_n\to K(\a)$ is injective, then
for every double coset $D\in T\backslash \Sy_n/G_\a$,
$$
  \mmaD(F(\b)) = 0
  \quad\iff\quad
  \b
  = [\phi(a_1),...,\phi(a_n)]^\sigma \>\>\text{for some $\sigma\in D$}.
$$
\end{proposition}

\begin{proof}
We have that
$\mmaD(F(\b)) = 0$ if and only if
$F(\b)=\phi(x)$ for some root $x$ of $\mmaD$ in $K(\a)$.
Such roots are $\eFa(\sigma)$ for some $\sigma\in D$, so
\beq
  \mmaD(F(\b)) = 0 &\iff& F(\b)=\phi(\eFa(\sigma)) &\text{for some $\sigma\in D$}\cr
       &\iff& F(\phi^{-1}(\b))=\eFa(\sigma)=F(\a^\sigma)\cr
       &\iff& \phi^{-1}(\b) = (\a^\sigma)^\tau = \a^{\tau\sigma}
              &\text{for some $\tau\in T$}\cr
       &\iff& \b = \phi(\a^{\sigma'}) = \phi(\a)^{\sigma'}
              &\text{for some $\sigma'\in D$}.\cr
\eeq
\end{proof}

\begin{theorem}
\label{cor1}
Let $F$ be a $G_\a$-invariant with
$\eFa\!\!: G_\a\backslash\,\Sy_n\!\to\! K(\a)$~\hbox{injective}.
If $F(\b)=F(\a)\in K$, then $a_i\mapsto b_i$ defines an isomorphism
$K(\a)\to K(\b)$.
\end{theorem}

\begin{proof}
Take $T=G_\a$ and $D$ the principal double coset $G_\a1G_\a$, and apply
the proposition. Since $\mmaD(X)=X-F(\a)$,
we have $\mmaD(F(\b)) = 0$, so $\b=\phi(\a)^\sigma$ for some
$\sigma\in G_\a$ and some isomorphism
$\phi: K(\a)\to K(\b)$.
Then $\phi\circ\sigma$ is the required isomorphism.
\end{proof}

\begin{remark}
\label{prob1sol}
This gives a solution to Problem \ref{problem1}:

Pick a $G_\a$-invariant $F$, e.g. using Remark \ref{invariantexists}.
Adjusting $f$ if necessary, we may assume
that $\eFa: T\,\backslash\,\Sy_n\to K(\a)$ is injective
(Lemma \ref{ugly}c).
In $L'$, keep permuting the roots of $f$ until $F(\b)$
becomes $F(\a)\in K$.
When this happens, $a_i\mapsto b_i$ defines an isomorphism
of the two splitting fields.

Note however, that in the worst case
we are evaluating a polynomial with $|G|$ terms on $|G\backslash \Sy_n/G|$
permutations. So the complexity is about $n!$ operations,
which is impractical for large $n$.
\end{remark}

\begin{example}[$\Di_{10}$-extensions]
Suppose $f(x)\in K[x]$ has degree 5, and $G_\a=\Gal(f/K)$
is the dihedral group $\Di_{10}$, generated by $(12345)$ and $(25)(34)$.
Take
$$
  F(x_1,...,x_5) = x_1 x_2 + x_2 x_3 + x_3 x_4 + x_4 x_5 + x_5 x_1.
$$
This is a $T$-invariant with $T=G_\a$: it is clearly invariant under $\Di_{10}$,
and on the other hand a permutation preserving $F$ is determined by
$x_1\mapsto x_i$, $x_2\mapsto x_{i\pm 1}$, so there are at most 10 choices.
In particular, $F(a_1,...,a_5)$ is invariant under the Galois group, and so
lies in $K$. Substituting the $a_i$ into $F$ in all possible orders
gives the values
$$
  \e\a F(\sigma^{-1}) =
  a_{\sigma(1)} a_{\sigma(2)} + a_{\sigma(2)} a_{\sigma(3)} +
  a_{\sigma(3)} a_{\sigma(4)} + a_{\sigma(4)} a_{\sigma(5)} +
  a_{\sigma(5)} a_{\sigma(1)}.
$$
Clearly each one occurs at least $10$ times for varying $\sigma\in \Sy_5$,
corresponding to the fact that
$\e\a F$ factors through $\Di_{10}\backslash \Sy_5$.
The assumption that the map \hbox{$\eFa: T\,\backslash\,\Sy_n\to K(\a)$} is injective
simply says that there are no more repetitions, and there are $120/10=12$
distinct values.

Suppose that this is indeed
the case, and let $b_1,...,b_5$ be the roots of $f$ in some other splitting
field. If we substitute the $b_i$ in $F$ in all possible orders~$\b^\sigma$,
we get again 12 values, one of which is $F(a_1,...,a_5)\in K$.
There are 10 isomorphisms $K(\a)\to K(\b)$ obtained from one another by
composing with Galois. They are determined by $\a\mapsto \b^\sigma$
for 10 permutations $\sigma\in \Sy_n$. Clearly, for each of these $\sigma$,
we have $F(\b^\sigma)=F(\a)$.
But, since every value is taken
exactly 10 times, we have the converse as well: if $F(\b^\sigma)=F(\a)$
for some $\sigma\in \Sy_n$, then $\a\mapsto\b^\sigma$ must define an
isomorphism of the splitting fields.
So to find an isomorphism, we only need to locate $F(\a)$ among the 12
values $F(\b^\sigma)$.

Note that the other values $F(\b^\sigma)$ are not in general $K$-rational,
so we cannot compare them with the values on $\a$. Their minimal
polynomials are the $\mm\a DF(X)$ for the 4 double cosets
$\Di_{10}\backslash \Sy_5/\Di_{10}$.
\end{example}

\section{Recognising conjugacy in Galois groups}
\label{sreccon}

In questions such as computing Frobenius elements in Galois groups
it is not necessary to compare the roots in two splitting fields.
It suffices to identify the conjugacy class of a specific Galois
automorphism:

\begin{problem}
\label{problem2}
Let $f(x)\in K[x]$ be a separable polynomial which splits completely
in $L\supset K$, and suppose we know $G=\Gal(f/K)$ as a permutation
group on the roots in $L$. If $L'$ is another field where $f$
splits completely and we are
given a permutation of the roots of $f$ in $L'$ which comes from some
Galois automorphism, find the conjugacy class of this automorphism in $G$.
\end{problem}

\begin{remark}
An isomorphism $\phi$ of the two splitting fields of $f$ induces an
isomorphism of Galois groups $G$ and $G'$.
We would like to identify
an element $\cB\in G'$ as an element $\cA\in G$. Note that $\cA$ depends
on the choice of~$\phi$. As any two isomorphisms differ by a Galois automorphism,
the conjugacy class $[\cA]$ is well-defined and this is what we are after.
\end{remark}

It is easy to see that a solution to Problem~\ref{problem1} answers
Problem~\ref{problem2} as well, so this is a weaker question.
However, we aim for a more practical solution (see Remark \ref{prob1sol}).
We may clearly restrict our attention to one cycle type in $\Sy_n$.
For convenience, throughout the section we we also fix a representative:

\begin{notation}
Fix an element $\xi\in \Sy_n$ and write $Z_\xi\< \Sy_n$ for its centraliser.
\end{notation}

\begin{definition}
\label{defcM}
Suppose $\Psi\in \Sy_n$ is conjugate to $\xi$, in other words they have
the same cycle type, say $\xi=\sigma_0 \Psi \sigma_0^{-1}$.
For a $T$-invariant $F$ and an ordering $\a$ of the roots of $f$,
define the polynomial
$$
  \cM F\a\Psi(X) = \prod_{\sigma\in (Z_\xi\cap T)\backslash Z_\xi\sigma_0}\mm\a\sigma F(X).
$$
It is well-defined by Lemma \ref{eFprop}(3).
Note that $Z_\xi\sigma_0$ is the set of all permutations that conjugate
$\Psi$ to $\xi$, in particular it is independent of the choice of $\sigma_0$.
\end{definition}

\begin{remark}
The situation we have in mind is that we have two sets of roots
$\a$ and $\b$ of $f$ in different splitting fields. So
there is an isomorphism
$\phi\!:\! K(\a)\!\to\!K(\b)$,
but we do not have it explicitly.
However, suppose~we~
know
that an automorphism $\cA\in\Gal(K(\a)/K)$
corresponds to $\cB\in\Gal(K(\b)/K)$ under $\phi$, and that
they permute the roots by
$$
  \cA(\a) = \a^\Psi, \quad  \cB(\b) = \b^\xi, \qquad \Psi,\xi\in \Sy_n.
$$
Then $\{\a^\sigma\}_{\sigma\in Z_\xi\sigma_0}$
is the set of all reorderings of $\a$ on which $\cA$ acts as $\xi$,
and $\cM F\a\Psi(X)$ is the smallest $K$-rational polynomial that has
$F(\a^\sigma)$ as roots for all such $\sigma$. But $\phi^{-1}(\b)$
must be one of these reorderings because $\cB$ acts on $\b$ as $\xi$.
The upshot is that $\cM F\a\Psi(X)$ has $F(\b)$ as a root,
and its construction does not require the knowledge of $\phi$.
In other words, if $\cM F\a\Psi(F(\b))\ne 0$, then we know that $\cA$ does
not correspond to $\cB$ under any isomorphism.
(In~\S\ref{sts} we will take $T=Z_\xi$ and turn this into an if and
only if statement.)
\end{remark}

\begin{lemma}
\label{lemFrobab}
Let $\phi: K(\a)\to K(\b)$ be an isomorphism of two splitting fields of $f$,
and define $\rho\in \Sy_n$ by $\b = \phi(\a^\rho)$. Then
$$
  \cM F\a{\rho^{-1}\Phi\rho} = \cM F\b\Phi.
$$
\end{lemma}

\begin{proof}
Write $\Psi=\rho^{-1}\Phi\rho$. Pick
$\sigma_\Phi$ with $\xi=\sigma_\Phi\Phi\sigma_\Phi^{-1}$, and
let $\sigma_\Psi=\sigma_\Phi\rho$, so that
$$
  \sigma_\Psi\Psi\sigma_\Psi^{-1}
    = \sigma_\Phi\rho\Psi\rho^{-1}\sigma_\Phi^{-1}
    = \sigma_\Phi\Phi\sigma_\Phi^{-1}
    = \xi.
$$
By definition,
$$
  \cM F\b\Phi \>\,= \!\!\!\prod_{\scriptscriptstyle\sigma\in (Z_\xi\cap T)\backslash Z_\xi\sigma_\Phi} \!\!\!\mm\b\sigma F,
    \qquad\quad
  \cM F\a\Psi \>\, = \!\!\!\prod_{\scriptscriptstyle\sigma\in (Z_\xi\cap T)\backslash Z_\xi\sigma_\Psi} \!\!\!\mm\a\sigma F.
$$
We claim that
$$
  \mm\a{s\sigma_\Psi}F = \mm\b{s\sigma_\Phi}F\qquad \text{for $s\in Z_\xi$}.
$$
First we show that they have the same degree.
Because $G_\b=\rho G_\a\rho^{-1}$ by the definition of $\rho$,
\beq
  \deg\mm\a{s\sigma_\Psi}F &=&
  |T\backslash Ts\sigma_\Psi G_\a| =
  |T\backslash Ts\sigma_\Psi G_\a\rho^{-1}| \cr
  &=& |T\backslash Ts\sigma_\Phi\rho G_\a\rho^{-1}| =
  |T\backslash Ts\sigma_\Phi G_\b|
  = \deg\mm\b{s\sigma_\Phi}F.
\eeq
Since both polynomials are powers of irreducible ones,
it now suffices to identify one of the roots:
\beq
  \e\a F(s\sigma_\Psi)
  &=& \e\a F(s\sigma_\Phi\rho)
  = F(\a^{s\sigma_\Phi\rho}))
  = F(\phi^{-1}(\b)^{s\sigma_\Phi}))\cr
  &=& F(\phi^{-1}(\b^{s\sigma_\Phi}))
  = \phi^{-1}(F(\b^{s\sigma_\Phi}))
  = \phi^{-1}(\e\b F(s\sigma_\Phi)).
\eeq
\end{proof}

\begin{corollary}
\label{cor4.18}
The map $\Psi\mapsto \cM F\a\Psi$ is constant on every conjugacy class
of $G_\a$ with cycle type $\xi$.
\end{corollary}

\begin{proof}
By the lemma above, $\cM F\a\Psi = \cM F\a{g\Psi g^{-1}}$ for $g\in G_\a$.
\end{proof}

We now have an approach to Problem \ref{problem2}:

\begin{proposition}
\label{proprec}
Let $\a, \b$ be orderings of the roots of $f$ in two different splitting
fields, and suppose $\Psi\in G_\a$ and $\Phi\in G_\b$ have cycle type $\xi$.
If the polynomials $\cM F\a\psi$ are distinct for $\psi$ in
different conjugacy classes of $G_\a$ of cycle type $\xi$, then
$$
  \vabove{\raise-5pt\hbox{\text{
    there is an isomorphism $K(\a)\to K(\b)$
  }}}{\raise5pt\hbox{\text{
    under which $\Psi$ corresponds to $\Phi$
  }}} \quad\Longleftrightarrow\quad
  \cM F\a\Psi=\cM F\b\Phi.
$$
If, moreover, the $\cM F\a\psi$ are pairwise coprime, then this occurs
precisely when $\cM F\a\Psi(F(\b^\sigma))=0$ for
some (any) $\sigma\in \Sy_n$ with $\xi=\sigma\Phi\sigma^{-1}$.
\end{proposition}

\begin{proof}
`$\Rightarrow$' is Lemma \ref{lemFrobab}. For `$\Leftarrow$',
pick any isomorphism $\phi: K(\a)\to K(\b)$.
The
polynomial $\cM F\b\Phi$ agrees with some $\cM F\a\psi$ by the lemma,
and $\Psi$ lies in the conjugacy class of $\psi$ by assumption.
Composing $\phi$ with an automorphism of $K(\a)/K$
(which corresponds to conjugating $\psi$) we obtain the required isomorphism.
\end{proof}

\begin{example}[Serre's trick \cite{Buh,Rob}]
\label{serre1}
Suppose $\vchar K\ne 2$,
$f\in K[x]$ has degree~$n$, and
$G_\a=\Gal(f/K)$ is the alternating
group $\Alt_n$. There is a particularly nice $T$-invariant with $T=\Alt_n$,
a `square root of the discriminant'
$$
  F(x_1,...,x_n) = \prod_{i<j} (x_i-x_j).
$$
The only double cosets $TxG_\a$ in $\Sy_n$ are $D=\Alt_n$ and
its complement $D'$ in~$\Sy_n$. Clearly
$\mm\a{D}F(X)=X-F(\a)$ and $\mm\a{D'}F(X)=X+F(\a)$, and $F(\a)^2=\Delta_f$
is the discriminant of $f$. So if $\b$ is the list of roots of $f$ in
some other splitting field, we find that
$$
  \vabove{\text{$a_i\mapsto b_i$ defines an }}
         {\text{isomorphism $K(\a)\to K(\b)$}}
  \quad\iff\quad
\prod_{i<j} (a_i-a_j)=\prod_{i<j} (b_i-b_j).
$$
This illustrates Theorem \ref{cor1} in the case of $\Alt_n$. To explain
Proposition \ref{proprec} in this setting,
suppose $\xi\in \Sy_n$ is a product of cycles of distinct odd degrees, so
that there are two conjugacy classes $[\Psi_1], [\Psi_2]$ in $G_\a=\Alt_n$
of~cycle type~$\xi$ (e.g. 5-cycles in $\Alt_5$).
Say $\sigma_1\Psi_1\sigma_1^{-1}\!=\!\xi\!=\!\sigma_2^{}\Psi^{}_2\sigma_2^{-1}$
with $\sigma_1\!\in\! \Alt_n$~and~$\sigma_2\!\notin\! \Alt_n$.
In this case $Z_\xi\subset \Alt_n=T$, so
\beq
  \cM F\a{\Psi_1}(X) = \mm\a{\sigma_1}F(X) = \mm\a {D\phantom{'}}F(X) = X-F(\a), \\[2pt]
  \cM F\a{\Psi_2}(X) = \mm\a{\sigma_2}F(X) = \mm\a {D'}F(X) = X+F(\a). \cr
\eeq
Suppose again that $\b$ is the list of roots of $f$ in
some other splitting field, and $\cB\in\Gal(K(\b)/K)$ is an automorphism
of cycle type $\xi$. Rearranging the $b_i$ if necessary, assume that $\cB$
acts on the $b_i$ as $\xi$, i.e. $\cB(\b)=\b^\xi$.
The statement of the proposition is that
$$
  \vabove{\text{$\cB$ comes from $[\Psi_1]$ under an}}
         {\text{isomorphism $K(\a)\to K(\b)$}}
  \quad\iff\quad
  \prod_{i<j} (a_i-a_j)=\prod_{i<j} (b_i-b_j),
$$
which is precisely Serre's trick.
The same invariant $F$ may sometimes be used in other subgroups of $\Sy_n$
to distinguish between the conjugacy classes of such cycle types.
(It determines whether the two classes are conjugate in $\Alt_n$ or not.)
\end{example}

\section{The directed edges invariant}
\label{sts}

As before, suppose $f(x)\in K[x]$ is separable
and $\a=[a_1,...,a_n]$ are its (ordered) roots in a splitting field.
We apply the results of \S\ref{sreccon} when $T=Z_\xi$, the centraliser of $\xi$.
This is particularly nice for two reasons: first, the polynomials
$\cM F\a\psi$ of Proposition \ref{proprec} are irreducible and distinct, and
second, it is easy to write down a $T$-invariant with just $n$ terms and
of degree~3 (compare the polynomials in Remark \ref{invariantexists}
and Example \ref{TSinv}).

\begin{proposition}
\label{propTS}
Let $\xi\!\in\!\Sy_n$ with centraliser $Z_\xi$. Suppose that $F$ is
a \hbox{$Z_\xi$-invariant}
such that $\eFa: Z_\xi\,\backslash\,\Sy_n\to K(\a)$ is injective.
Let $\Psi, \Psi'\in G_\a$ be two elements of cycle type $\xi$. Then
\begin{enumerate}
\item $\cM F\a{\Psi}$ is irreducible, and equals $\mm\a{\sigma}F$
      for any $\sigma\in\Sy_n$ with $\xi=\sigma \Psi \sigma^{-1}$.
\item $\cM F\a{\Psi}$ has degree $|[\Psi]|$.
\item $\cM F\a{\Psi}=\cM F\a{\Psi'}$ if and only if $\Psi$ and $\Psi'$
are conjugate in $G_\a$.
\end{enumerate}
\end{proposition}

\begin{proof}
For brevity, write $Z\!=\!Z_\xi$.
Pick $\sigma, \sigma'\!\in\!\Sy_n$ with
$\sigma \Psi \sigma^{-1}\!=\!\xi\!=\!\sigma' \Psi (\sigma')^{-1}$.

\noindent
(1) By definition,
$$
  \cM F\a{\Psi} \>=\>
    \prod_{\tau\in (Z\cap Z)\backslash Z\sigma}\!\!\! \mm\a{\tau}F
  \>\>\>=\> \mm\a{\sigma}F.
$$
It is irreducible by the assumed injectivity of $\eFa$ (see Remark \ref{reminj}).

\noindent
(2) By definition,
\beq
  \deg\mm\a{\sigma}F &=& \displaystyle
    |Z\backslash Z\sigma G_\a| = \frac{|Z\sigma G_\a|}{|Z|}
      = \frac{|\sigma^{-1} Z\sigma G_\a|}{|Z|}\cr
    &=& \displaystyle \frac{|G_\a|}{|G_\a\cap \sigma^{-1} Z\sigma|}
    = \frac{|G_\a|}{|\Cent_{G_\a}(\Psi)|}=|[\Psi]|.
\eeq
(3) If $\Psi$ and $\Psi'$ are conjugate, then $\cM F\a{\Psi}=\cM F\a{\Psi'}$
by Corollary \ref{cor4.18}. Conversely, suppose that
$\cM F\a{\Psi}=\cM F\a{\Psi'}$.
Since $\eFa$ is injective,
$Z\sigma G_\a=Z\sigma' G_\a$, so
$\sigma'=s\sigma g$ for some $s\in Z$ and $g\in G_\a$. Then
$$
  \Psi' = (\sigma')^{-1} \xi \sigma'
         = g^{-1}\sigma^{-1} s^{-1} \xi s\sigma g
         = g^{-1}\sigma^{-1} \xi \sigma g
         = g^{-1}\Psi g,
$$
so $[\Psi']=[\Psi]$.
\end{proof}

\begin{example}[The directed edges invariant]
\label{TSinv}
Let $\xi\in \Sy_n$ and fix a polynomial $h\in K[x]$ of degree at least 2.
Define
$$
  F(x_1,...,x_n)=\sum_{j=1}^n h(x_j)\, x_{\xi(j)}.
$$
It can be visualised as the directed edges in a graph that
define the action by $\xi$.
For instance, for $\xi=(1234)(56)\in \Sy_6$ and $h(x)=x^2$,
$$
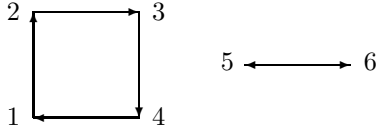
\begin{picture}(180,80)
  \put(20,27){\smaller[1]1}
  \put(75,27){\smaller[1]4}
  \put(20,67){\smaller[1]2}
  \put(75,67){\smaller[1]3}
  \put(101,48){\smaller[1]5}
  \put(155,48){\smaller[1]6}
  \put(30,30){\vector(0,1){40}}
  \put(30,70){\vector(1,0){40}}
  \put(70,70){\vector(0,-1){40}}
  \put(70,30){\vector(-1,0){40}}
  \put(130,50){\vector(1,0){20}}
  \put(130,50){\vector(-1,0){20}}
  \put(-20,5){$F\>=\>x_1^2x_2\!+\!x_2^2x_3\!+\!x_3^2x_4\!+\!x_4^2x_1\>\>+\>\>x_5^2x_6\!+\!x_6^2x_5$}
\end{picture}
$$
It is clearly a $Z_\xi$-invariant.
\end{example}

\begin{definition}
\label{gammaC}
Fix $h(x)\in K[x]$. For each conjugacy class $C$ in $G_\a$ define
\label{defgamma}
$$
  \m_C(X) = \prod_{\sigma\in C} (X-\sum_{j=1}^n h(a_j) \sigma(a_j)).
$$
\end{definition}

\begin{lemma}
\label{mjST}
Let $F$ be as in Example \ref{TSinv}. Then for every $\Psi\in G_\a$,
$$
  \cM F\a{\Psi}(X) = \m_{[\Psi]}(X).
$$
\end{lemma}

\begin{proof}
Pick $\sigma\in\Sy_n$ with $\sigma \Psi \sigma^{-1}=\xi$.
First, suppose $\tau\in [\Psi]$ and $u_\tau\in\Sy_n$ satisfies
$u_\tau^{-1}\xi u_\tau=\tau$. Then
\beq
  \eFa(u_\tau) &=&\displaystyle F(\a^{u_\tau}) = \sum_i h(a_{u_\tau^{-1}(i)})a_{u_\tau^{-1}(\xi(i))}\cr
    &=&\displaystyle \sum_j h(a_j) a_{u_\tau^{-1}\xi u_\tau(j)}
    = \sum_j h(a_j) \tau(a_j).
\eeq
On the other hand, note that for $t\in Z_\xi$ and $g\in G_\a$,
$$
  (t\sigma g)^{-1}\xi(t\sigma g)= g^{-1}\sigma^{-1}t^{-1} \xi t\sigma g =
  g^{-1}\sigma^{-1} \xi \sigma g = g^{-1} \Psi  g.
$$
So for $\tau = g^{-1} \Psi  g \in [\Psi]$,
$$
  \{u_\tau\in S_n\>|\>u_\tau^{-1}\xi u_\tau=\tau\} = Z_\xi \sigma g,
$$
because the left-hand side is clearly some right coset of $Z_\xi$.
This equality gives a correspondence between $[\Psi]$ and
$Z_\xi\backslash Z_\xi\sigma G_\a$. So
\beq
  \cM F\a{\Psi}(X) &=& \displaystyle \mm\a{\sigma }F(X) =
  \prod_{u\in (Z_\xi\backslash Z_\xi\sigma G_\a)} (X-\eFa(u)) \cr
   &=& \displaystyle
  \prod_{\tau\in [\Psi] } (X-\sum_{j=1}^n h(a) \tau(a_j)) = \m_{[\Psi]}(X),
\eeq
as claimed.
\end{proof}

\begin{corollary}
\label{corfrobred}
Let $\a, \b$ be orderings of the roots of $f$ in two different splitting
fields, and let $\Psi\in G_\a$ and $\Phi\in G_\b$.
If the $\m_C(X)$ are pairwise coprime for
different conjugacy classes of $G_\a$, then
$$
  \vabove{\raise-5pt\hbox{\text{
    \smaller[1]
    there is an isomorphism $K(\a)\to K(\b)$
  }}}{\raise5pt\hbox{\text{
    \smaller[1]
    under which $\Psi$ corresponds~to~$\Phi$,
  }}} \>\>\Longleftrightarrow\>\>
  \m_{[\Psi]}\Bigl({\textstyle\sum_j h(b_j) \Phi(b_j)}\Bigr)=0.
$$
The condition that the $\m_C$ are coprime is satisfied for
$h(x)$ in a Zariski dense open set
in the space of all polynomials of degree at most $n-1$.
\begin{proof}
The equivalence follows from Proposition \ref{proprec} and the lemma above.
For the last assertion apply Lemma \ref{ugly2}.
\end{proof}

\end{corollary}

\section{Frobenius elements}
\label{sfrob}

Now suppose $K$ is a global field. We turn to our initial problem
of computing Frobenius elements in Galois groups.
We use the following remarkable property
of the directed edges invariant:

\begin{proposition}
\label{frobred}
Let $f(x)\in \F_q[x]$ be a polynomial with roots
$a_1,...,a_n\in\bar\F_q$ counted with multiplicity,
and let $\phi=\Frob_q\in\Gal(\bar\F_q/\F_q)$. For every poly\-nomial
$h(x)\in \F_q[x]$,
$$
  \sum_{j=1}^n h(a_j) \phi(a_j) = \tr_{A/\F_q}(h(X)X^q),
$$
where $X$ is the class of $x$ in the algebra $A=\F_q[x]/f$.
\end{proposition}
\noindent
This is an immediate consequence of the lemma below (with $H(x)=h(x)x^q$).
\begin{lemma}
\label{lemfrobred}
Let $k$ be a field and $f(x)\in k[x]$ a polynomial with roots
$a_1,...,a_n\in\bar k$ counted with multiplicity.
Then for every $H(x)\in k[x]$,
$$
  \sum_{j=1}^n H(a_j) = \tr_{A/k}(H(X)),
$$
where $X$ is the class of $x$ in $A=k[x]/f$.
\end{lemma}

\begin{proof}
Consider $X$ as a linear map $A\to A$, $Y\mapsto XY$. Its minimal polynomial
is $f$, since $f(X)=0$ but no linear combination of $1,X,...,X^{n-1}$
is zero. So the generalised eigenvalues of $X$ are exactly the $a_i$,
and those of $H(X)$ are therefore $H(a_i)$ (look at the Jordan normal form
of $X$ over $\bar k$).
The result follows.
\end{proof}

\begin{theorem}[Generalised Euler's criterion]
\label{frobmain}
Let $K$ be a global field and $f(x)\in K[x]$ a separable
polynomial with roots $a_1,...,a_n$ in $\bar K$ and Galois group $G$.
Fix $h(x)\in K[x]$ and for each conjugacy class $C$ of $G$ set
$$
  \m_C(X) = \prod_{\sigma\in C} (X- \sum_{j=1}^n h(a_j)\sigma(a_j)).
$$
\begin{itemize}
\item[(a)]
The polynomials $\m_C(X)$ have coefficients in $K$.
\item[(b)]
Let $\p$ be a prime of $K$ with residue field $\F_q$, and $C$ a conjugacy
class of $G$. If $\p$ does not divide
the denominators of the coefficients of $f$ and $h$,
the leading coefficient of $f$ and
the resultants $\Res(\m_C,\m_{C'})$ for $C'\ne C$,
then the coefficients of $\m_C(X)$ are integral at $\p$ and
$$
  \Frob_\p\in C \quad\iff\quad
  \m_C\Bigl(\tr_{\frac{\F_q[x]}{f(x)}/\F_q}(h(x)x^q)\Bigr)=0\mod\p.
$$
\item[(c)]
For all $h(x)$ in some Zariski dense open set
in the space of polynomials of degree at most $n-1$,
we have $\Res(\m_C,\m_{C'})\ne 0$ for every pair of conjugacy classes
$C\ne C'$.
\end{itemize}
\end{theorem}

\begin{proof}
(a) This follows from Lemma \ref{mjST}, Definition \ref{defcM}
and Remark \ref{reminj}.

(b) $\m_C(X)$ is clearly integral at the required primes.

`$\Rightarrow$': if $\Frob_\p\in C$ then
$\sum_{j=1}^n h(a_j) \Frob_\p(a_j)$ is a root of $\m_C(X)$ by
the definition of $\m_C$, and it reduces mod $\p$ to
${\tr_{\frac{\F_q[x]}{f(x)}/\F_q}(h(x)x^q)}$ by Proposition~\ref{frobred}.

`$\Leftarrow$': the polynomial $\m_C(X)$ is distinguished from the
others by any one of its root mod $\p$
by the assumption that $\p\nmid\Res(\m_C,\m_{C'})$ for $C\ne C'$.

(c) Apply Lemma \ref{ugly2}.
\end{proof}

\begin{remark}[Choice of $h$]
\label{remres}
If the resultants $\Res(\m_C,\m_{C'})$ are non-zero,
Theorem \ref{frobmain}b describes the Frobenius element for all but finitely
many primes~$\p$.
If one of the resultants vanishes, equivalently $\m_C$ has a common
factor with some $\m_{C'}$, the statement does not apply to $C$ for any $\p$.
However, this is rare and easily avoided by
choosing a different $h$; most choices will work by Theorem \ref{frobmain}c.

Alternatively, for any fixed $h$ with $1<\deg h<n$ it is possible to
replace $f$ by another polynomial $\tilde f$ of degree $n$ with
the same splitting field so that the resulting $\m_C$ are coprime. To see this,
consider
$$
  \gamma_C(X) = \prod_{\sigma\in C} (X-\sum_{j=1}^n h(x_j) x_{\sigma(j)}),
$$
and note that they are coprime as polynomials in $X$ over $K(x_1,...,x_n)$.
Now apply Lemma \ref{ugly}b to $F_1=\prod_{C\ne C'}\Res(\gamma_C,\gamma_{C'})$
and $F_2=0$. We obtain a Zariski dense open set of polynomials $B(t)$
of degree at most $n-1$ for which $\tilde f=\prod_j(x-B(a_j))$ works.
\end{remark}

\endcomment

\begin{remark}[Euler's criterion]
The classical criterion
\smash{$a^{\frac{p-1}2}\equiv(\frac ap)\mod p$}
says that \smash{$a^{\frac{p-1}2}=\pm 1$}
determines whether $x^2-a$ has a root modulo $p$.
Similarly, to see whether $x^3-a$ has a root modulo $p\equiv 1\mod 3$
one checks whether \smash{$a^{\frac{p-1}3}$} is 1 or another third root
of unity in $\F_p^\times$, etc.

One can reformulate this as a matrix statement: take a $2\times 2$
matrix $M$ with minimal polynomial $x^2-a$ (respectively $3\times 3$ and
$x^3-a$). Then $M^{p-1}$ is the scalar matrix with $a^{\frac{p-1}2}$
(respectively $a^{\frac{p-1}3}$) on the diagonal, so its trace determines
whether the polynomial has a root in $\F_p$; e.g. for $x^3-a$ the
distinction
is whether $\frac13\tr M^{p-1}$ is $1$ or a root of $x^2+x+1$.

Theorem \ref{frobmain} generalises this to arbitrary polynomials over
global fields. Observe that for a polynomial
$$
  f(x)=x^n+c_{n-1}x^{n-1}+\ldots+c_0
$$
the trace in the theorem can be interpreted as a trace of a matrix, e.g.
\begingroup
\baselineskip 8pt
\def\arraystretch{0.8}
\def\arraycolsep{3pt}
$$
  \tr_{\frac{\F_q[x]}{f(x)}/\F_q}(x^d) = \tr
  \begin{pmatrix}
   & &&-c_0  \cr
  1&&&-c_1  \cr
   &\raise3pt\hbox{$\ddots$}&&\raise3pt\hbox{$\vdots$}\cr
   & &1\!&-c_{n-1} \cr
  \end{pmatrix}^{\rlap{\raise-3pt\hbox{$\!\!d$}}}\mod q.
$$
\endgroup
Therefore (a minor modification of) the trace $\Tr M^{q-1}$ for a
matrix $M$ with minimal polynomial $f$ determines the splitting behaviour
of $f$ mod $\p$ and the conjugacy class of Frobenius, in the same way
as above. See also \S\ref{sintro} and~\S\ref{sexanab}.
\end{remark}

\begin{remark}[Ramified primes]
\label{rampr}
The condition that $\p$ does not divide any resultant
$\Res(\m_C,\m_{C'})$ excludes
all primes that ramify in the splitting field of $f$ over $K$. Indeed,
if $\sigma\ne 1$ is an element of inertia at $\q$ for some $\q|\p$, it
is easy to see that $\m_{[1]}$ and $\m_{[\sigma]}$ have a common root mod $\p$.
\end{remark}

\begin{remark}[Extending to all $\p$]
\label{FrobQp}
In order to deal with the primes dividing the resultants,
we may work over the completion $K_\p$ instead
of the residue field $\F_q$.
Compute the splitting field $L/K_\p$ of $f$ and the roots $b_1,..,b_n$.
Choose a lift $\Psi$ of the Frobenius element in $\Gal(L/K_\p)$ and evaluate
$$
  \sum_{j=1}^n h(b_j) \Psi(b_j).
$$
This number is now a root of precisely one of the $\m_C$, and this $C$
is the conjugacy class of the chosen Frobenius lift $\Psi$.
(See Corollary \ref{corfrobred}.)
\end{remark}

\begin{remark}[Artin L-functions]
\label{artinl}
Suppose $L/K$ is a Galois extension of number fields with Galois group $G$,
represented as a splitting field of some polynomial
$f(x)\in K[x]$. Recall that a complex representation $\rho$ of $G$
is called an {\em Artin representation\/}. It has an L-series
defined by the Euler product over all primes of $K$,
$$
  L(\rho,s) = \prod_\p \frac{1}{P_\p(q^{-s})}.
$$
Here $q$ is the size of the residue field at $\p$ and
$$
  P_\p(T) = \det (1-\Frob_\p T \>|\> \rho^{I_\p})
$$
is the inverse characteristic polynomial of Frobenius on the subspace of
$\rho$ fixed by the inertia group $I_\p$ at $\p$.

Theorem \ref{frobmain} and Remark \ref{FrobQp} allow us to explicitly
compute the coefficients of such L-series.
For the unramified primes, they recover
the conjugacy class of $\Frob_\p$ in $G$, which determines the local
polynomial $P_\p(T)$. For the ramified primes, it suffices to find
the restriction of $\rho$ to the local Galois group $G_\p$ at $\p$
with respect to an embedding $G_\p\injects G$ as a decomposition group.
Assuming we can find $G_\p$, Remark \ref{FrobQp} enables us to
identify the conjugacy class in $G$ of any element of $G_\p$, under this
embedding. This is sufficient to compute the character of $\rho$ on $G_\p$,
and thus also $\rho^{I_\p}$ and $P_\p(T)$.
Note that we have {\em not} actually found the decomposition group at $\p$
as a subgroup of $G$,
which appears to be a harder problem.

This algorithm to compute Frobenius elements and L-series of Artin
representations has now been implemented in Magma \cite{Magma}.
\end{remark}

\begin{remark}[Complexity]
From the complexity point of view, the computation of Frobenius elements
for `good' primes has two steps:

One is the initial precomputation of the polynomials $\m_C$,
each of which takes $O(n|C|)$ operations in some field containing the
$a_j$ (e.g. $\C$ or $\bar\Q_p$). This needs to be done for all conjugacy
classes that are not determined by their cycle type.

The second step deals with a specific prime $\p$ of $K$ with residue
field~$\F_q$. We determine the cycle type of $\Frob_\p$ by computing
$\gcd(f,x^{q^j}-x)$ for $j\le n/2$, which takes $O(n\log q)$
multiplications of $n\times n$ matrices over $\F_q$.
Then we evaluate the trace $\tr(h(x)x^q)$
with another $O(n+\log q)$ matrix multiplications.
Finally, we substitute the trace into all $\m_C$ corresponding to the
cycle type of $\Frob_\p$, which is $O(d)$ coefficient reductions and
multiplications in $\F_q$, where $d$ is the number of elements
in $G$ of this cycle type.

Here is as an illustration for polynomials of degree at most 11.
There are 474 transitive groups $G$ on at most $11$ points, for each of which
we took a polynomial $f\in\Q[x]$ with $\Gal f=G$ as a permutation group on
the roots. (We used the database in Magma \cite{Magma} V2.16.)
For each $G$ we computed $\Frob_p$ for all $p<100000$ with $p\nmid\Delta_f$,
using Serre's trick (Example \ref{serre1}) and the above algorithm.
Together with the Galois group computation and the precomputation of the
$\m_C$ this took under 15 seconds on a 2GHz Pentium notebook for each $G$,
with only four exceptions that took longer:
$G=\Alt_5^2\rtimes \Cy_2$,
$\Alt_5^2\rtimes \Cy_2^2$,
$\Alt_5^2\rtimes \Cy_4$ and $\Ma_{11}$.
\end{remark}

\begin{remark}[Additional symmetries]
\label{remsym}
Suppose all conjugacy classes of elements
of some order $o$ and a fixed cycle type are closed
under the power maps $g\mapsto g^k$ for $k$ in some non-trivial subgroup
$H\subset(\Z/o\Z)^\times$ (for instance they are self-inverse, like in
dihedral groups). Then one may replace $\m_C(X)$ in Theorem \ref{frobmain} by
$$
  \prod_\sigma\Bigl(X-\sum_{j=1}^n h(a_j) (\mathop{\Sigma}\limits_{\scriptscriptstyle k\in H}\sigma^k(a_j))\Bigr),
$$
taking the product over some representatives for $C$ modulo
the action of~$H$, and modifying the trace accordingly.
In practice, this speeds up the computation of the $\m_C$, as their degree
drops by a factor of $|H|$.
\end{remark}

\vskip 3mm

\section{Examples: abelian groups}
\label{sexaab}

\vskip 2mm

If the Galois group is abelian, its conjugacy classes are of size 1,
and all the $\m_C$ of Theorem \ref{frobmain} are linear,
$\m_C(X)=X-r_C$ with $r_C\in K$.
For a good choice of $h(x)$ and all but finitely many primes $\p$,
the trace $\tr(h(x)x^q)$ agrees with exactly one of the $r_C$ modulo $\p$,
which then determines the conjugacy class of $\Frob_\p$.

In the examples below, $\zeta_n$ denotes a primitive $n$th root of unity.

\begin{example}
\label{exqi}
Let $K=\Q(i)$ and
$$
  f(x) = x^4 + 2x^3 + (3+3i)x^2 + 4ix -1+i.
$$
Its complex roots are
$a_1\!=\!-0.31795\!-\!0.57510i$, \
$a_2\!=\!0.50870\!-\!1.1289i$,
\hbox{$a_3\!=\!-1.4682\!+\!1.8471i$} and
$a_4\!=\!-0.72255\!-\!0.14308i$
to 5 decimal places. The~splitting field $L$ is a $\Cy_4$-extension of $\Q(i)$, non-Galois
over $\Q$, and the Galois group of $L/K$ is
$\langle(1234)\rangle\<\Sy_4$. Take $h(x)=x^2$. An elementary computation
gives
\beq
  \m_{[\id]}      &=& X-(10+6i), &&&
  \m_{[(1234)]}   &=& X-(4+4i), \cr
  \m_{[(13)(24)]} &=& X-(-2+2i), &&&
  \m_{[(1432)]}   &=& X+8. \cr
\eeq
For a prime $\p\ne (1+i), (2-i), (3)$ (the primes dividing $r_C-r_{C'}$
for $C\ne C'$) with residue field $\F_q$, we deduce that the Frobenius at $\p$
is determined~by
$$
\def\SM#1{#1}
\begin{array}{|c|ccccc|}
\hline
  \SM{
  \tr_{\frac{\F_q[x]}{f(x)}/\F_q}({x^{q+2}})\equiv} &&
    10\!+\!6i & 4\!+\!4i & -2\!+\!2i & -8 \vphantom{\int^X}\cr
\hline
  \hfill\Frob_\p\hfill= && \id & (1234) & (13)(24) & (1432) \cr
\hline
\end{array}
$$
\end{example}

\begin{example}[Kummer extensions]
\label{exkummer}
Suppose $\zeta=\zeta_n\in K$ and $L=K(\sqrt[n]s)$ is a Kummer extension
of degree $n$. It is abelian with Galois group $\Cy_n$ whose elements
are determined by
$$
  \sigma_i: \sqrt[n]s \quad\longmapsto\quad \zeta^i \sqrt[n]s, \qquad
  i=1,\ldots,n.
$$
Take $f(x)=x^n-s$ and $h(x)=x^{n-1}$. 
Then
\beq
  \m_{[\sigma_i]}(X) =
     \displaystyle X - \sum_{j=1}^n h(\zeta^j\sqrt[n]s)\sigma_i(\zeta^j\sqrt[n]s) =
     X - ns\cdot\zeta^i.
\eeq
For a prime $\p$ of $K$ with residue field $\F_q$,
because \hbox{$n\>|\>q\!-\!\!1$}, we have
$$
  \tr_{\frac{\F_q[x]}{f(x)}/\F_q}(h(x)x^q) =
  \tr_{\frac{\F_q[x]}{x^n\!-\!s}/\F_q}(x^{q+n-1}) =
  \tr_{\frac{\F_q[x]}{x^n\!-\!s}/\F_q}(s^{\frac{q-1}n+1}) =
  ns\cdot s^{\frac{q-1}n}.
$$
So Theorem \ref{frobmain} says that for $\p\nmid ns$,
$$
  \Frob_\p=\sigma_i \quad\iff\quad s^{\frac{q-1}n}\equiv\zeta^i\mod\p,
$$
which is the classical criterion for Kummer extensions.
\end{example}

\begin{example}[$\Q(\zeta_p)/\Q$]
Let $\zeta=\zeta_p$ for some prime $p>2$, and take
$$
  K=\Q, \quad L=\Q(\zeta), \quad f(x)=x^{p-1}+\ldots+x+1.
$$
Thus $\Gal(L/K)\iso(\Z/p\Z)^\times$, with elements
$\sigma_i: \zeta \mapsto \zeta^i$ for $i=1,\ldots,p-1$.
For $h(x)=x^2$ we have $\m_{[\sigma_i]}(X)=X-r_i$ with $r_i\in\Q$ given by
$$
  r_i = \sum_{j=1}^{p-1} (\zeta^j)^2 \sigma_i(\zeta^j)
                  = \sum_{j=1}^{p-1} \zeta^{j(2+i)} =
         \bigleftchoice {-1}{i\ne p-2,}{p-1}{i=p-2.}
$$
For a prime $q$ of $\Q$,
\beq
  \tr_{\frac{\F_q[x]}{f(x)}/\F_q}(h(x)x^q)
  &=& \tr_{\frac{\F_q[x]}{f(x)}/\F_q}(x^{q+2})
  &\equiv& \tr_{\frac{\Z[x]}{f(x)}/\Z}(x^{q+2})\mod q \\[2pt]
  &\equiv& \tr_{F/\Q}(\zeta^{q+2})
  &\equiv&\leftchoice {-1}{p\nmid q+2}{p-1}{p|q+2}\mod q.
\eeq
Hence Theorem \ref{frobmain}b shows that for all $q\ne p$,
$$
  \Frob_q = \sigma_{p-2} \quad\liff\quad q\equiv -2\mod p.
$$
The same computation with $h(x)=x^{p-k}$ for varying $k$ yields
the classical criterion
$$
  \Frob_q = \sigma_{k} \quad\liff\quad q\equiv k\mod p.
$$
Note that none of these $h(x)$ work for all conjugacy classes
simultaneously, because the $\m_{[\sigma_j]}$ are not coprime.
This tends to happen when
the roots of $f$ are `too nice' and $h(x)$ is `too simple'. By Lemma
\ref{ugly2}, most $h$ do work. In our example, a general polynomial
$$
  h(x) = \lambda_{1} x^{p-1} + \ldots + \lambda_{p-1} x + \lambda_p
$$
has
$$
  \m_{[\sigma_i]}(X) = X + h(1) - p \lambda_i,
$$
and these are distinct if and only if $\lambda_1,\ldots,\lambda_{p-1}$ are.
The primes to which the theorem then applies are those not dividing
$p\prod(\lambda_i-\lambda_j)$ in this case.
\end{example}

\begin{example}[Cyclotomic extensions]
In general, suppose $L=K(\zeta_n)$ is some cyclotomic extension,
and $f(x)$ is the minimal polynomial of $\zeta_n$ over~$K$.
As in the previous example, $G=\Gal(L/K)\injects (\Z/n\Z)^\times$, and we
write $\sigma_i$ for the automorphism with $\sigma_i(\zeta_n)=\zeta_n^i$.
We do the same computation as above: for $h(x)=x^k$ and $\p$
a prime of $K$ with residue field $\F_q$,
$$
  \m_{[\sigma_i]}(X)
     = X \!-\! \sum_{g\in G} g(\zeta_n)^k \sigma_i(g(\zeta_n))
     = X \!-\! \sum_{g\in G} g(\zeta_n)^{k+i}
     = X \!-\! \Tr_{L/K}(\zeta_n^{k+i})
$$
and
$$
  \tr_{\frac{\F_q[x]}{f(x)}/\F_q}(x^{k+q})
   \equiv \tr_{L/K}(\zeta_n^{k+q})\mod \p.
$$
Because $\tr_{L/K}(\zeta_n^j)$ is $|G|$ precisely when $n|j$,
the polynomial $\m_{[\sigma_{n-k}]}$ differs from all the other $\m_{[\sigma_j]}$'s,
and we find that
$$
  \Frob_\p = \sigma_{n-k} \quad\liff\quad q\equiv n-k\mod n
$$
for almost all $\p$.
(One may improve `almost all' to `all $\p\nmid n$' by taking several~$h$.)
\end{example}

\begin{remark}
The fact that we obtained a simple formula for Frobenius elements for
cyclotomic and Kummer extensions relied on the existence of a universal
expression for the trace $\tr(h(x)x^q)\mod \p$. It follows from class
field theory that there are such formulae in all abelian extensions.

For instance, consider Example \ref{exqi} of a $\Cy_4$-extension of $K=\Q(i)$ from
the point of view of class field theory.
There the conductor of $L/K$ is $N=(1+i)^4(2-i)=8-4i$, and
the group $(\O_K/N)^\times$ is $\Cy_4\times \Cy_4\times \Cy_2$, with generators
$i$, $7$ and $3-2i$ respectively.
For a prime $\p=(\alpha)\subset\Z[i]$
not dividing $N$, if $\alpha\equiv i^a 7^b (3-2i)^c \mod N$, then
$\Frob_\p=(1234)^b$.

Now compare this with the description of Frobenius in Example \ref{exqi}.
Writing $\F_q=\Z[i]/\p$ and $\tr$ for $\tr_{\frac{\F_q[x]}{f(x)}/\F_q}$,
we get 4 congruences for the traces,
\beq
  \p=(\alpha),\>\> \alpha\equiv i^a 7^0 (3-2i)^c \mod N &\iff&
  \tr({x^{q+2}})\equiv 10+6i\mod \p  \\[2pt]
  \p=(\alpha),\>\> \alpha\equiv i^a 7^1 (3-2i)^c \mod N &\iff&
  \tr({x^{q+2}})\equiv 4+4i\mod \p  \\[2pt]
  \p=(\alpha),\>\> \alpha\equiv i^a 7^2 (3-2i)^c \mod N &\iff&
  \tr({x^{q+2}})\equiv -2+2i\mod \p  \\[2pt]
  \p=(\alpha),\>\> \alpha\equiv i^a 7^3 (3-2i)^c \mod N &\iff&
  \tr({x^{q+2}})\equiv -8\mod \p
\eeq
for $\p\ne (1+i),(2-i),(3)$.

Note that if one had a way to prove these congruences directly, one would have a proof
of Artin reciprocity in the extension $L/K$.
\end{remark}

\section{Examples: non-abelian groups}
\label{sexanab}

We continue with examples to Theorem \ref{frobmain}.
When $G$ is non-abelian, the only difference is that the $\m_C$
are no longer linear.

\begin{example}
Let $K=\Q$ and $f(x)=x^3-2$. It has Galois group $\Sy_3$ and roots
$a_1=\sqrt[3]{2}, a_2=\zeta\sqrt[3]{2}$ and $a_3=\zeta^2\sqrt[3]{2}$, where
$\zeta$ is a primitive cube root of unity. Take $h(x)=x^2/6$
(the factor $1/6$ is only chosen for convenience)
and compute the polynomials
$\m_C$ for the three conjugacy classes:

\begingroup
\vskip-3mm
\smaller[2]
\beq
\m_{[\id]}&=X-\tfrac16(a_1^2a_1+a_2^2a_2+a_3^2a_3)\\[2pt]
  &=X-1\\[2pt]
\m_{[(12)]}\!\!\!&=
  (X\!-\!\tfrac16(a_1^2a_2\!+\!a_2^2a_1\!+\!a_3^3))
  (X\!-\!\tfrac16(a_1^2a_3\!+\!a_2^3\!+\!a_3^2a_1))
  (X\!-\!\tfrac16(a_1^3\!+\!a_2^2a_3\!+\!a_3^2a_2))\\[2pt]
  &=(X-\tfrac13(\zeta+\zeta^2+1))(X-\tfrac13(\zeta^2+1+\zeta))(X-\tfrac13(1+\zeta+\zeta^2))\\[2pt]
  &=X^3\\[2pt]
\m_{[(123)]}\!\!\!\!\!&=
  (X-\tfrac16(a_1^2a_2+a_2^2a_3+a_3^2a_1))
  (X-\tfrac16(a_1^2a_3+a_2^2a_1+a_3^2a_2))\\[2pt]
  &=(X-\tfrac13(\zeta+\zeta+\zeta))(X-\tfrac13(\zeta^2+\zeta^2+\zeta^2))=(X-\zeta)(X-\zeta^2)\\[2pt]
  &=X^2+X+1.
\eeq
\endgroup
On the other hand, for a rational prime $q=3m+k$ with $k=1$ or 2,
$$
  \tr_{\frac{\F_q[x]}{x^3-2}/\F_q}(\tfrac16 x^{q+2}) =
  \tr(\tfrac16 2^{m+1}x^{k-1})
  =\leftchoice{2^m}{k=1}{0}{k=2}
  =\leftchoice{2^{\frac{q-1}3}}{q\equiv1\mod 3}{0}{q\equiv2\mod 3}.
$$
The conclusion of Theorem \ref{frobmain} is that, as expected, for $q\ne 2,3$,
\beq
q\equiv1\mod 3,\>\> 2\in(\F_q)^{\times3} & \implies & \Frob_q=\id,\cr
q\equiv1\mod 3,\>\> 2\notin(\F_q)^{\times3} & \implies & \Frob_q\in[(123)],\cr
q\equiv2\mod 3& \implies & \Frob_q\in [(12)].\cr
\eeq
Clearly, an identical computation goes through for $f(x)=x^3-c$ (with $h(x)=x^2/3c$)
over any global field $K$ with $\zeta\not\subset K$.
\end{example}

We can also take a general cubic polynomial and obtain an analogue of
Euler's criterion for its factorisation modulo primes:

\begin{theorem}
\label{gencubus}
Let $f(x)=x^3+bx+c$ be a separable cubic polynomial over a global field $K$,
and $\p$ a prime of $K$ with residue field $\F_q$.
Write
$$
  T = \tr_{\frac{\F_q[x]}{f(x)}/\F_q}(x^{q+1}) =
      \tr \begin{pmatrix}0&0&-c\cr1&0&-b\cr0&1&0\end{pmatrix}^{\!\!q+1}\mod \p.
$$
If $\p$ does not divide $3b(4b^3+27c^2)$ and the denominators of $b$ and $c$,
then
$$
\begin{array}{lllllll}
  T&\equiv&-2b\!\!\!&\mod \p&\iff&\text{$f(x)$ has 3 roots mod $\p$},\cr
  T&\equiv& b &\mod \p&\iff&\text{$f(x)$ is irreducible mod $\p$},\cr
  \multicolumn{4}{l}{\text{$T$ is a root of $x^3\!-\!3b^2x\!-\!2b^3\!-\!27c^2$}}
  &\iff&\text{$f(x)$ has 1 root mod $\p$}.\cr
\end{array}
$$
\end{theorem}

\begin{proof}
We compute the polynomials $\m_C$ for $G=\Sy_3$, $h(x)=x$ by
expressing their
coefficients in terms the elementary symmetric functions
$a_1\!+\!a_2\!+\!a_3\!=\!0$, $a_1a_2\!+\!a_2a_3\!+\!a_3a_1=b$ and $a_1a_2a_3=-c$:

\begingroup
\vskip-2mm
\smaller[1]
\beq
\m_{[\id]} &=\!\!& X-(a_1^2\!+\!a_2^2\!+\!a_3^2)=X-(a_1\!+\!a_2\!+\!a_3)^2\!+\!2(a_1a_2\!+\!a_1a_3\!+\!a_2a_3) \cr
           &=\!\!& X + 2b\\[1pt]
\m_{[(12)]}&=\!\! & (X\!-\!(a_1a_2\!+\!a_2a_1\!+\!a_3^2))(X\!-\!(a_1a_2\!+\!a_2a_1\!+\!a_3^2))(X\!-\!(a_1a_2\!+\!a_2a_1\!+\!a_3^2))\cr
           &=\!\!&X^3-3b^2X-2b^3\!-\!27c^2\\[1pt]
\m_{[(123)]}\!\!&=\!\! & (X-(a_1a_2\!+\!a_2a_3\!+\!a_3a_1))(X-(a_1a_3\!+\!a_2a_1\!+\!a_3a_2))\cr
           &=\!\!& (X-b)^2.
\eeq
\endgroup
The least common multiple of their pairwise resultants is $3b(4b^3+27c^2)$,
which completes the proof by Theorem \ref{frobmain}.
\end{proof}

An identical computation can be done for polynomials of higher degree,
as long as one has the patience to work out the coefficients of the~
$\m_C$'s.
Here is the corresponding result for quartics:

\begin{theorem}
Let $f(x)=x^4+bx^2+cx+d$ be a separable quartic polynomial over $K$,
and $\p$ a prime of $K$ with residue field $\F_q$.
Then the value $\tr_{\frac{\F_q[x]}{f(x)}/\F_q}(x^{q+1})$ is a root of
one of the polynomials

\begingroup
\vskip-1mm
\smaller[2]
\beq
\m_{[\id]}&= &
X+2b   \cr
\m_{[(12)(34)]} &=&
X^3-2bX^2-16dX+32bd-8c^2 \cr
\m_{[(12)]} &=&
X^6+4bX^5+(2b^2+8d)X^4+(-12b^3+48bd-26c^2)X^3\cr
&&-(23b^4-120b^2d+108bc^2+112d^2)X^2\cr
&&-(16b^5-128b^3d+138b^2c^2+256bd^2+216c^2d)X\cr
&&-4b^6+48b^4d-56b^3c^2-192b^2d^2-288bc^2d-27c^4+256d^3 \cr
\m_{[(123)]} &=&
X^4+(-2b^2+8d)X^2-8c^2X+b^4-8b^2d+8bc^2+16d^2 \cr
\m_{[(1234)]} &=&
X^3-2bX^2+(b^2-4d)X+c^2.\cr
\eeq
\endgroup
If $\p$ does not divide the denominators of $b,c$ and $d$ and the
pairwise resultants of the $\m_c$, then this determines the degrees
in the factorisation of $f$ mod $\p$: they are the cycle lengths of
the permutation in the index of $\m$.
\end{theorem}

A theorem of Brumer (see \cite{JLY} Thm. 2.3.5) states that any Galois
extension $L/K$ with Galois group $G=\Di_{10}$ is a splitting field of
$$
  f_{a,b}(x) = x^5 + (a\!-\!3)x^4 + (b\!-\!a\!+\!3)x^3
               + (a^2\!-\!a\!-\!1\!-\!2b)x^2 + bx + a
$$
for some $a,b\in K$. Using a similar argument to $G=\Sy_3$ and $\Sy_4$,
we find

\begin{theorem}
Suppose $L/K$ is the splitting field of $f_{a,b}(x)$ as above, with
$G=\Gal(L/K)\iso \Di_{10}$.
If $\p$ a prime of $K$ with residue field $\F_q$, not dividing
$3a-b+1$ and the denominators of $a$ and $b$ and such that $f\mod\p$
is irreducible, then $\tr_{\frac{\F_q[x]}{f(x)}/\F_q}(x^{q+1})$
is either $-2a+b+1$ or $a+2$ modulo $\p$.
This determines which of the two conjugacy classes
of 5-cycles contains $\Frob_\p$.
\end{theorem}


\begin{remark}
In this setting, if $\Frob_\p$ is not a 5-cycle, it is either
the identity or an element of order 2. In the former case,
\smash{$\tr_{\frac{\F_q[x]}{f(x)}/\F_q}(x^{q+1})$} is $a^2-4a-2b+3\mod \p$; in
the latter it is a root of

\begingroup
\vskip-2mm
\smaller[3]
\beq
  \m_{\scriptscriptstyle[(23)(45)]}=\!\!\!\!\!\!&
  X^5
\!-\!(a\!-\!3)^2X^4
\!+\!(31\!-\!2a^3\!+\!4b\!-\!3b^2\!+\!a^2(11\!+\!2b)\!-\!2a(21\!+\!2b))X^3
\cr
&\!+\!(12a^3(3\!+\!2b)\!-\!a^2(137\!+\!44b)\!+\!
a(114\!+\!6b\!-\!28b^2)\!-\!51\!+\!7a^4\!-\!4a^5\!-\!20b\!+\!14b^2\!-\!2b^3)X^2
\cr
&\!+\!(40\!+\!16a^5\!-\!8a^6\!+\!32b\!-\!17b^2\!-\!4b^3\!+\!a^4(58\!+\!42b)\!+\!
a^2(182\!+\!18b\!-\!52b^2)\cr&\!+\!4a^3(\!-\!49\!-\!21b\!+\!b^2)
\!-\!2a(65\!+\!13b\!-\!17b^2\!+\!6b^3))X
\cr
&\!+\!8a^6\!-\!4a^7\!+\!4a^5(7\!+\!5b)\!-\!4a^4(32\!+\!17b)\!+\!
2a^3(123\!+\!85b\!+\!4b^2)\cr& \!-\!a^2(245\!+\!218b\!+\!24b^2)\!-\!
2a(\!-\!30\!-\!6b\!+\!51b^2\!+\!22b^3)\!+\!2(\!-\!6\!-\!8b\!+\!3b^2\!+\!b^3\!-\!4b^4).
\eeq
\endgroup
\end{remark}

\begin{example}
Here is another example, to illustrate what the $\m_C$ look
like in general.
Take $K=\Q$ and $L=\Q(E[3])$,
the 3-torsion field of
the elliptic curve $E: y^2+y=x^3-x^2$.
Then $\Gal(L/K)\iso\GL_2(\F_3)$, and $L$ is the splitting field~of
$$
  f(x) = x^8-9x^7+18x^6+33x^5-93x^4-15x^3-23x^2-36x-27.
$$
The $\m_C$ for $h(x)=x^2$ are
\begingroup
\smaller[3]
\beq
\m_{\scriptscriptstyle [\id]} &=\!\!\!\!&
X\!-\!144\cr
\m_{\scriptscriptstyle [(1 3)(2 4)(5 6)(7 8)]}\!\!\!\!\! &=\!\!\!\!&
X\!-\!3\cr
\m_{\scriptscriptstyle [(2 4)(5 7)(6 8)]} &=\!\!\!\!&
X^{12}\!-\!699X^{11}\!+\!204666X^{10}\!-\!32922129X^9\!+\!3212225793X^8\!-\!196600821903X^7\!+\!\cr
  &=\!\!\!\!&7340079612456X^6\!-\!145234777501584X^5\!+\!566948224573848X^4\!+\!\cr
  &=\!\!\!\!&26747700562448082X^3\!-\!187604198442957555X^2\!-\!2946247136394353892X\!-\!\cr
  &=\!\!\!\!&24290099658154516203\cr
\m_{\scriptscriptstyle [(1 4 8)(2 7 3)]} &=\!\!\!\!&
X^8\!-\!546X^7\!+\!120102X^6\!-\!14088342X^5\!+\!989228043X^4\!-\!43566817716X^3\!+\!\cr
&=\!\!\!\!&1248800990265X^2\!-\!21583664066961X\!+\!167939769912993\cr
\m_{\scriptscriptstyle [(1 4 3 2)(5 7 6 8)]} &=\!\!\!\!&
X^6\!-\!258X^5\!+\!26448X^4\!-\!1344378X^3\!+\!34859664X^2\!-\!445164021X\!+\!2926293624\cr
\m_{\scriptscriptstyle [(1 7 4 3 8 2)(5 6)]} &=\!\!\!\!&
X^8\!-\!264X^7\!+\!29292X^6\!-\!1698042X^5\!+\!51288993X^4\!-\!654852960X^3\!+\!\cr
&=\!\!\!\!&3360584547X^2\!-\!277935306777X\!+\!7299371089503\cr
\m_{\scriptscriptstyle [(1 5 4 7 3 6 2 8)]} &=\!\!\!\!&
X^6\!-\!258X^5\!+\!26250X^4\!-\!1336755X^3\!+\!35700471X^2\!-\!477465444X\!+\!2707751520\cr
\m_{\scriptscriptstyle [(1 6 4 8 3 5 2 7)]} &=\!\!\!\!&
X^6\!-\!258X^5\!+\!28230X^4\!-\!1674048X^3\!+\!57362760X^2\!-\!1097286921X\!+\!9616023198\cr
\eeq
\endgroup
\end{example}

\begin{example}
\label{pgsp43}
As an indication to the kind of Artin L-series that may be
numerically computed, we give an example with a big Galois group over~$\Q$.
We take $G=\PGSp(4,\F_3)$ of order 51840, realised through the
Galois action on the 3-torsion of the Jacobian of a genus 2 curve, and
evaluate the Artin L-series of an irreducible 6-dimensional representation
of $G$.

Specifically, $G$ is the unique double cover of
the simple group $\Sp(4,\F_3)/\F_3^\times$
in $\PGL(4,\F_3)=\GL(4,\F_3)/\F_3^\times$.
To obtain it as a Galois group, take the hyperelliptic curve
$$
  \cC/\Q:\>\> y^2 - (x^2+1)y = x^5-x^4+x^3-x^2.
$$
Consider the field $\Q(J[3])$ obtained by adjoining to $\Q$ the
coordinates of the 3-torsion points of its Jacobian $J/\Q$.
Then $\Gal(\Q(J[3])/\Q)$ is $\GSp(4,\F_3)$.
The group we want is $G=\GSp(4,\F_3)/\{\pm 1\}$, and it can be obtained
from the Galois action on the 40 lines through the origin in $J[3]$.
Specifically, if $(P)\!+\!(Q)\!-\!2(O)\in J[3]$ is a non-zero point
with $P=(x_P,y_P), Q=(x_Q,y_Q)$,
the minimal polynomial $f$ of $x_Px_Q$ over $\Q$ has Galois group $G$;
$$
  f=x^{40} \!+\! 27x^{39} \!+\! 39x^{38} \!-\! 61x^{37} \!+\! \ldots
  + 2259x^{3} \!+\! 3471x^{2} \!+\! 1057x \!+\! 69.
$$
In its action on the roots of $f$, the group has several conjugacy classes
of the same cycle type, and the largest $\m_C$ that we need has degree 2160
(using Remark \ref{remsym}).

The group has two irreducible 6-dimensional representations, $\rho$ and $\rho'$
(whose trace on elements of order 10 in $G$ is $+1$ and $-1$ respectively).
The curve $\cC$ has good reduction outside 2 and 3, so $L/\Q$ is unramified
at all primes $p\ne 2,3$. The conductor of $\rho$ is $2^{10}3^{17}$
and we used our machinery to compute
the local polynomials for the Artin L-series $L(\rho,s)$
for primes up to 410203. Using Magma \cite{Magma}, we then evaluate
$$
  L(\rho, 1) \approx 1.852529796, \qquad L(\rho, 2) \approx 1.119877506
$$
to 10 digits precision; this computation relies implicitly on the validity
of Artin's conjecture for $\rho$.
The total time to compute $f$, $\Gal(f/\Q)$, the $\m_C$,
the L-series and the L-values
was 7 hours on a Sun Ultra 24 workstation.

\end{example}

\section{Appendix: Two lemmas on Zariski density}
\label{sugly}

\begin{lemma}
\label{ugly}
Suppose $K$ is an infinite field, $f\in K[t]$ is a separable polynomial of degree $n$
and $a_1,...,a_n$ are its roots in some splitting field $L$.
\begin{itemize}
\item[(a)]
If $F,G\in K[x_1,...,x_n]$ take the same values on
\beq
  x_1= \beta_0 + \beta_1 a_1 + \ldots + \beta_{n-1} a_1^{n-1} \cr
  \hskip 4mm
  \cdots\cr
  x_n= \beta_0 + \beta_1 a_n + \ldots + \beta_{n-1} a_n^{n-1} \cr
\eeq
for all $[\beta_1,...,\beta_n]\in K^n$, then $F=G$.
\item[(b)] Suppose $F_1,...,F_d\in K[x_1,...,x_n]$ are distinct. There exists
a polynomial $B(t)=\beta_0+\ldots+\beta_{n-1}t^{n-1}\in K[t]$ such that
$B(a_1),...,B(a_n)$ generate $L$ and the $F_i$ take distinct values on $[B(a_1),...,B(a_n)]$.
The set of such $B$ is Zariski dense in $K\oplus Kt\oplus\cdots\oplus Kt^{n-1}$.

\item[(c)] Let $F$ be a $T$-invariant for some $T\<\Sy_n$.
There is a Zariski dense open set of polynomials
$B(t)\in K\oplus Kt\oplus\cdots\oplus Kt^{n-1}$ for which
$\a'=[B(a_1),...,B(a_n)]$ generate $L$ and
$\e{\a'}{F}: T\,\backslash\,\Sy_n\to L$ is injective.
\end{itemize}
\end{lemma}

\begin{proof}
(a) Let $U=K(t_1,...,t_n)$. As a first step, we observe that $K^n$
is Zariski dense in $\A^n_U=U^n$: this is clear for $n=1$ as $K$ is infinite;
generally, if $K^n$ were not Zariski dense, it would be contained in
a (not necessarily irreducible) hypersurface of some degree $d$, so it
would contain at most $d$ hyperplanes. But, by induction, it contains
all $\{r\}\times U^{n-1}$ for all $r\in K$, which gives a contradiction.

Therefore, as $F$ and $G$ are continuous in the Zariski topology,
they agree on all of $U^n$, i.e. on all the above combinations
with $[\beta_1,...,\beta_n]\in U^n$. Now solve the
system of equation $\sum_{j=0}^{n-1}a_i^j\beta_j=t_j$ for $\beta_1,...,\beta_n$.
(This is possible because $a_i\ne a_k$ for $i\ne k$, so the
Vandermonde matrix is invertible.)
Using this solution
we find that $F(t_1,...,t_n)=G(t_1,...,t_n)$, so $F=G$ as polynomials.

(b) Put $F(x_1,...,x_n)=\prod_{i<j}(x_i-x_j)(F_i-F_j)$ and $G=0$ and apply~(a).
This gives a polynomial $B(t)=\beta_0+\ldots+\beta_{n-1}t^{n-1}\in K[t]$
which clearly satisfies the `distinct values' condition. Furthermore,
$B(a_i)\ne B(a_j)$ guarantees the `generate $L$' condition as well:
the Galois action permutes the $B(a_i)$ in the same way as the $a_i$,
so the Galois group has the same order.
Finally, consider $F(B(a_1),...,B(a_n))$ as a polynomial in
$\beta_0,...,\beta_{n-1}$.
Its zero set is Zariski closed in $\A^n$ and we proved that its
complement is non-empty. This proves the last claim.

(c) Apply (b) to the set of polynomials
$\{F^\sigma\}_{\sigma\in T\backslash \Sy_n}$, using that, by definition,
$\e{\a'}{F}(\sigma^{-1}) = F((\a')^{\sigma^{-1}}) = F^\sigma (\a')$.
\end{proof}

\begin{lemma}
\label{ugly2}
Suppose $K$ is an infinite field, $f\in K[t]$ is a separable polynomial
of degree $n$ and $a_1,...,a_n$ are its roots in some splitting field $L$.
Then on a Zariski dense open set of polynomials $h(x)$ in
$K\oplus Kx\oplus\ldots\oplus Kx^{n-1}\iso {\mathbb A}_K^n$, the values
$$
  v_h(\sigma) =  \sum_{j=1}^n h(a_j)\sigma(a_j),
    \qquad \sigma\in G=\Gal(L/K)
$$
are distinct.
\end{lemma}

\begin{proof}
For any $\sigma\in G$, the map $E_\sigma: h\mapsto v_h(\sigma)$
is $K$-linear $K^n\to L$. So $E_\sigma$ agrees with $E_\tau$
on a $K$-linear subspace for every $\sigma, \tau\in G$.
If none of these subspaces is all of $K^n$, then the complement of their
union is the desired set (non-empty since $K$ is infinite).
It remains to prove that $E_\sigma\ne E_\tau$ for~$\sigma\ne\tau$.

Suppose $E_\sigma=E_\tau: K^n\to L$. Then their extensions by linearity
to maps $L^n\to L$ agree as well. In other words, $v_h(\sigma)=v_h(\tau)$
for all $h$ in $L\oplus Lx\oplus\ldots\oplus Lx^{n-1}$. In particular,
taking
$$
  h(x) = \prod_{j\ne i} (x-a_j)
$$
we get that $\sigma(a_i)=\tau(a_i)$. As this holds for all $i$,
it follows that $\sigma=\tau$.
\end{proof}

\end{document}